\documentclass[12pt,leqno]{article}  
\usepackage{hyperref}
\usepackage{amsmath,amstext,amsthm,amssymb,amsrefs}   
\numberwithin{equation}{section}

\allowdisplaybreaks
\swapnumbers
\theoremstyle{plain}

\def\seq{\lesssim }

\def\bthm#1.#2 #3\ethm{
\begin{\ifx#1ttheorem\fi%
\ifx#1llemma\fi%
\ifx#1ccorollary\fi%
\ifx#1pproposition\fi%
\ifx#1ddefinition\fi}
\label{#1.#2}  
#3 \end{\ifx#1ttheorem\fi%
\ifx#1llemma\fi%
\ifx#1ccorollary\fi%
\ifx#1pproposition\fi%
\ifx#1ddefinition\fi}}

\def\t#1/{theorem~\ref{t#1}}   \def\Tcap#1/{Theorem~\ref{t#1}} 
\def\c#1/{corollary~\ref{c#1}}   \def\Ccap#1/{Corollary~\ref{c#1}} 
\def\l#1/{lemma~\ref{l#1}}        \def\Lcap#1/{Lemma~\ref{l#1}}  
\def\s#1/{section~\ref{s#1}}      
\def\e#1/{(\ref{e#1})}
\def\d#1/{definition~\ref{d#1}}
\def\f#1/{figure~\ref{f#1}}

\def\Label #1 {\label{#1}}

\def\norm#1.#2.{\lVert#1\rVert_{#2}}
\def\Norm#1.#2.{\bigl\lVert#1\bigr\rVert_{#2}}
\def\NOrm#1.#2.{\Bigl\lVert#1\Bigr\rVert_{#2}}
\def\NORm#1.#2.{\biggl\lVert#1\biggr\rVert_{#2}}
\def\NORM#1.#2.{\Biggl\lVert#1\Biggr\rVert_{#2}}

\def\ip#1,#2.{\langle #1,#2\rangle}
\def\Ip#1,#2.{\bigl\langle#1,#2\bigr\rangle}
\def\IP#1,#2.{\Bigl\langle#1,#2\Bigr\rangle}

\def\abs#1{\lvert#1\rvert}

\newcommand{\zb}{\ensuremath{\beta}}
\newcommand{\zc}{\ensuremath{\psi}}

\newcommand{\ze}{\ensuremath{\epsilon}}
\newcommand{\zve}{\ensuremath{\varepsilon}}
\newcommand{\zvf}{\ensuremath{\varphi}}

\newcommand{\zG}{\ensuremath{\Gamma}}
\newcommand{\zh}{\ensuremath{\eta}}
\newcommand{\zI}{\ensuremath{\infty}}

\newcommand{\zo}{\ensuremath{\theta}}

\newcommand{\zx}{\ensuremath{\xi}}

\newcommand{\zp}{\ensuremath{\pi}}
\newcommand{\zq}{\ensuremath{\chi}}
\newcommand{\zz}{\ensuremath{\zeta}}

\def\z#1#2{\ifcase#1 {{\mathcal {#2}}}  
\or {{\tilde{#2}}}                    
\or { {\boldsymbol{#2}}}                  
\or{{\widetilde{#2}}}                   
\or {{\acute{#2}}}			
\or \ensuremath{\ifx#2R\ZR\otimes\ZR\else\fi
		\ifx#2T\ZT\otimes\ZT\else\fi
		\ifx#2C\ZC_+\otimes\ZC_+\else\fi
		\ifx#2Z\ZZ\otimes\ZZ\else\fi
		\ifx#2D\ZD\otimes\ZD\else\fi
		\ifx#2N\ZN\otimes\ZN\else\fi}			
\or {{\bar{#2}}}
\or {\dot{#2}}
\or {\overline{#2}}
\or {\underline{#2}}\fi}

\def\zref #1.#2/{%
\ifx#1e(\ref{e.#2})\fi%
\ifx#1t{Theorem~\ref{t.#2}}\fi%
\ifx#1l{Lemma~\ref{l.#2}}\fi%
\ifx#1p{Proposition~\ref{p.#2}}\fi%
\ifx#1s{Section~\ref{s.#2}}\fi%
\ifx#1c{Corollary~\ref{c.#2}}\fi%
}

\def\ZR{\ensuremath{\mathbb R}}
\def\ZZ{\ensuremath{\mathbb Z}}
\def\ZD{\ensuremath{\mathbb D}}
\def\ZN{\ensuremath{\mathbb N}}
\def\ZT{\ensuremath{\mathbb T}}

\def\ZP{\ensuremath{\mathbb P}}
\def\ZC{\ensuremath{\mathbb C}}

\def\mid{\,:\,}

\def\md#1#2\emd{\ifx0#1
\begin{equation*} #2 \end{equation*}\fi  
\ifx1#1\begin{equation}#2\end{equation}\fi   
\ifx2#1\begin{align*}#2\end{align*}\fi   
\ifx3#1\begin{align}#2\end{align}\fi    
\ifx4#1\begin{gather*}#2\end{gather*}\fi  
\ifx5#1\begin{gather}#2\end{gather}\fi   
\ifx6#1\begin{multline*}#2\end{multline*}\fi  
\ifx7#1\begin{multline}#2\end{mutline}\fi  
}

\def\mo#1{\text{\rm#1MO}}
 
 \def\ind#1{ {\mathbf 1}_{#1}} 

\begin{document}

\title{Remarks on Product $\text{VMO}$ }
 \author{Michael T. Lacey\footnote{This work has been
  supported by an NSF grant.} 
\\ { Georgia Institute of Technology} \and 
Erin Terwilleger\footnote{Research supported in part by an NSF  VIGRE grant to the Georgia Institute of Technology.} \\
{University of Connecticut}
\and 
Brett D. Wick\footnote{Research supported in part by an NSF VIGRE grant to Brown University.}  \\ Brown University 
}

\date{}

\maketitle

\abstract{Well known results related to the compactness of Hankel operators of one complex 
variable are extended to little Hankel operators of two complex variables.    Critical to these 
considerations is the result of 
 Ferguson and Lacey \cite{sarahlacey} characterizing the boundedness of the little Hankel operators 
 in terms of the product \mo B of S.-Y.~Chang and R.~Fefferman \cites{cf1,cf2}.}

\parskip=12pt 

\section{Introduction}

We prove necessary and sufficient conditions for the compactness of little Hankel operators of two complex variables. 
In the one complex variable case, results of this type are sometimes referred to as Hartman's theorem. 
Central to this are the Hardy space and $\mo B$ of two complex variables. Formally, the easiest way to 
phrase these results is for the Hardy space $H^1(\z5T)$ and its dual space $\mo B(\z5T)$.   Definitions are 
postponed until the next section.

 $L^2(\z5T)$ is the direct sum of 
\md0
L^2(\z5T)=\oplus_{\zve\in\{\pm,\pm\}} H^2_{\zve}(\z5D)
\emd
in which $H^2_{\pm,\pm}(\z5D)$ is the space of square integrable functions with (anti) holomorphic extensions in 
each variable separately.  Let $\ZP_{\pm,\pm}$ be the corresponding projection of $L^2(\z5T)$ onto $H^2_{\pm,\pm}(\z5D)$.  

The Hankel operators
of interest to us are operators from $H^2_{+,+}(\z5D)$ to $H^2_{-,-}(\z5D)$ given by 
$h_\zvf:=\ZP_{-,-}M_{ \zvf}$ in which $M_\zvf$ denotes the operator of pointwise multiplication by $\zvf$.  
The following theorem extends Nehari's Theorem \cite{nehari} to two complex variables, and is  essentially a restatement of 
the main result of S.~Ferguson and the first author \cite{sarahlacey}.   We indicate a modification of the classical proof, which 
relies in an essential way on the results of \cite{sarahlacey}.

\bthm t.nehari  The Hankel operator $h_\zvf$ is bounded iff there is a function $\zc\in L^\zI(\z5T)$ for which 
$\ZP_{-,-}\zvf=\ZP_{-,-}\zc$, and we have the equivalence  
\md3 \label{e.nehari}
\norm h_\zvf ..\approx{}& \inf\{ \norm \zc.\zI.\mid \ZP_{-,-}\zvf=\ZP_{-,-}\zc\}
\\{}\approx{}& \norm  \ZP_{-,-}\zvf.\mo B(\z5D).  .  \label{e.nehari-bmo}
\emd
\ethm

The $\mo B$ space is the dual to real
$H^1(\z5D)$, as identified by S.-Y.~Chang and R.~Fefferman \cite{cf1,cf2}.  
We have the following refinement of this theorem.

\bthm c.compact  $h_\zvf$ is compact iff 
$\ZP_{-,-}\zvf$ is in the   closure of $C(\z5T)$ with respect to the $\mo B$ topology.
\ethm 

 In view of the classical result of Sarason \cite{sarason}, we call this last space $\mo V(\z5D)$.  This space has an equivalent characterization in terms of 
Carleson measures.  
In the circumstance in which the symbol is assumed to be bounded, we can say 
a little more.   Let $\z0L^p_{\pm,\pm}(\z5T)$ be the space of  functions $b\in L^p(\z5T)$ such that $\ZP_{\pm,\pm}b=0$. 

\bthm  t.0.1  
Let $\zvf\in L^{\infty}(\z5T)$.  Then the
following are equivalent
\begin{itemize}
\item[(i)]$h_\zvf$ is compact.\\
\item[(ii)] $\zvf\in\z0L^{\infty}_{-,-}(\z5T)+C(\z5T)$.\\
\item[(iii)] there exists a $g\in C(\z5T)$ such that $h_\zvf=h_g$.
\end{itemize}
\ethm

This theorem is a consequence of a finer fact about the 
essential norm of a little Hankel operator.  Take the essential norm to be  
\md0
\norm h_\zvf .e.:=\inf\{ \norm h_\zvf -K..\mid K\,:\,H_{+,+}(\z5T)\to{}H_{-,-}(\z5T)\text{ ${}$ is compact}\}.
\emd
Observe that $\norm h_\zvf.e.=0$ iff $h_\zvf$ is compact.
\bthm t.0.2
Let $\zvf\in L^{\infty}(\z5T)$.  Then
$$
\norm{h_{\zvf}}.e.\approx\textnormal{dist}_{L^{\infty}}(\zvf,\z0L^{\infty}_{-,-}+C).
$$
\ethm

These results have different, equivalent formulations in terms of Hankel matrices, or Hankel operators on 
$H^2(\z5C)$.   In addition, it is of interest to state a result in the equivalent language of commutators. Namely  for a function 
$\zvf\in \mo B(\z5R)$  define
\md1\label{e.cb}
C_\zvf:=[[M_\zvf,H_1],\, H_2]
\emd
in which $H_j$ denotes the Hilbert transform computed in the coordinate $j$. 
The Hilbert transform can be taken on the circle or the real line.  At this point, we take it to be defined on the real line.
Let us define 
\md0
\mo V(\z5R):=\textup{clos}_{\mo B }{C_0^\zI(\z5R)}
\emd
where $C_0^\zI$ denotes the space of smooth compactly supported functions. We will return to the Carleson measure characterization of membership in $\mo V$ later.

\bthm  t.many  
We have ${\mo V (\z5R)}^*=H^1(\z5R)$.  In addition, $C_\zvf$ is compact iff $\zvf\in\mo V(\z5R)$. 
\ethm

The next section contains background material for this paper.   Following that, the corollaries and theorems  related to compact operators 
are given in sections three and four.  The final section discusses the Carleson measure characterization of $\mo V$, and the duality 
statement $\mo V^*=H^1$.

\section{The Hardy Spaces of Two Complex Variables}

In speaking of Hardy spaces, one should take care to specify whether the functions are analytic, or not.   The analytic
Hardy spaces $H^p(\z5D)$ 
consists of functions $F\mid \z5D\to\ZC$ such that $F$ is holomorphic in each variable seperately, and 
\md0
\norm F.p.^p=\sup_{0<r_j<1}\int_{\z5T}\abs{ F(r_1e^{2\zp i\zo_1},r_2e^{2\zp i \zo_2})}^p\; d\zo_1\,d\zo_2. 
\emd
In the case that $1\le{}p\le\zI$, the boundry values of $F$ exist almost everywhere.  And the Fourier transform of that 
function is supported on the postive orthant of $\z5Z$.   In speaking of the analytic Hardy spaces, and their duals, we 
will use the notation $H^p(\z5D)$, $H^p(\z5C)$.

The (real) Hardy space $H^1(\z5R)$ consists of real valued functions $f$ on $\ZR^2$ for which 
\md0
\norm f.H^1(\z5R).:=\sum_{A_1,A_2\in\{I,H_1,H_2\}}\norm A_1A_2f.1.<\zI.
\emd
Here $I$ is the identity operator and $H_j$ is the Hilbert transform computed in the $j$th coordinate.  For $f\in H^1(\z5R)$, there is a biholomorphic extension $F(z_1,z_2)$ to $\z5C$ such that 
\md0
\lim_{y_1,y_2\downarrow0}\operatorname{Re} F(x_1+iy_1,x_2+iy_2)=f(x_1,x_2)\qquad \text{a.e.}
\emd
There are several equivalent definitions of this Hardy space in terms of maximal, square, and area functions, all formulated in terms of a product setting.
In speaking of the real Hardy spaces, we will use the notation $H^1(\z5R)$ or $H^1(\z5T)$. 

The dual space $\mo B(\z5R)$ was identified by S.-Y.~Chang and R.~Fefferman.  Their characterization is notable,
 as the structure of the allied Carleson measures is far more complicated than in a one parameter setting. This space has two known intrinsic characterizations.  One is that $\mo B(\z5R)$ is the dual to $H^1(\z5R)$, and the 
 second is that the $\mo B$ norm is comparable to the $L^2$ norm of the commutator $C_b$, which is one 
 formulation of the main result of Ferguson and Lacey \cite{sarahlacey}.   
 
 There is another definition in terms of wavelets and Carleson measures which, though no longer intrinsic in nature, is very useful.   We let $\z0D$ denote the set of  dyadic intervals on $\ZR$.  Given a rectangle $R=R_1 \times R_2\in\z0D\times\z0D$  define translation and dilation invariant operators by 
\md4
T_y f(y):=f(x-y),\qquad y\in\ZR^2,
\\
D^p_{R_1\times R_2} f(x_1,x_2):=\frac1{(\abs{R_1}\abs{R_2})^{1/p}} f\Bigl( \frac{x_1}{\abs{R_1}}, \frac{x_2}{\abs{R_2}}\Bigr)
,\qquad 0<p<\zI.
\emd  
Note that the second condition preserves $L^p$ norm and depends upon the scale but not location of the rectangle $R_1\times R_2$.

Given a function $w(x_1,x_2)=\prod_1^2v(x_j)$, we set 
\md0
w_R=T_{c(R)}D^2_R  w,\qquad \text{$c(R)={}$ the center of $R$}.
\emd
Our assumptions are  that $v$ is a  bounded, piecewise continuous, rapidly decreasing,  mean zero function, and that $\{w_R\mid R\in\z0D\times\z0D\}$ is an $L^2$ normalized orthogonal basis for $L^2(\ZR^2)$. 

Then, it is a theorem of Chang and Fefferman \cite{cf1,cf2} that we have 
\md1\label{e.CM}
\norm f . \mo B (\z5R).\approx{} \sup_U \Bigl[\abs {U}^{-1}\sum_{R\subset U} \abs{\ip f,w_R. }^2\Bigr]^{1/2}.
\emd
What is essential in this definition is that the supremum be formed over all subsets $U$ of the 
plane with finite measure.

To define analytic $\mo B(\ZC_+\otimes\ZC_+)$, one can use the same definition, provided one restricts attention to 
the jointly analytic projections of the wavelets. That is, the functions $w_R$ are replaced by $v_R:=\ZP_{+,+}w_R$, 
and then a definition of analytic $\mo B$ is just \zref e.CM/ with the $w_R$ replaced by $v_R$.

By $A\lesssim{}B$ we mean that there is an absolute constant $K$ so that $A\lesssim{}KB$.  By $A\approx{}B$ we mean $A\lesssim{}B$ and 
$B\lesssim{}A$.

\section{The Hankel Operators on $H^2(\z5D)$} 

\begin{proof}[Proof of \zref t.nehari/]

If it is the case that $\zc\in L^\zI(\z5T)$ exists with $\ZP_{-,-}\zvf=\ZP_{-,-}\zc$, then clearly we 
can estimate 
\md0
\norm h_\zvf..\le{} \norm \ZP_{-,-}{\zc} .2.\le{}\norm \zc.\zI. .
\emd
It is also then the case that $ \norm \ZP_{-,-} \zc .BMO(\z5D).\lesssim\norm \zc.\zI.$. 

\medskip

In the converse direction, we adopt a classical method of proof but use in an essential way the results 
of \cite{sarahlacey}.
We show that there is a $\zc\in L^\zI(\z5T)$ with  $\ZP_{-,-}\zvf=\ZP_{-,-}\zc$, and  
$\norm \zc.\zI.\lesssim{}\norm h_\zvf..$. We do so by defining a linear functional on  $H^1(\z5D)$ 
with norm less than or equal to a constant times $\norm h_\zvf..$.  For a pair of functions $f,g\in H^2(\z5D)$, set 
\md0
L(fg)=\int (h_\zvf f)g\; dx=\int (P_{-,-}\zvf) fg \; dx.
\emd
It is important to observe that this definition does not depend upon the order in which $f$ and $g$ are given to us.
And in addition, 
$\abs{ L(fg)}\le{}\norm h_\zvf ..\norm f.H^2.\norm g.H^2.$. Therefore, this 
definition of $L$ extends to the injective tensor product 
$H^2(\z5D)\widehat\otimes H^2(\z5D)$ which has the norm 
\md0 
\norm h. H^2\widehat\otimes H^2.:=\inf\Bigl\{ \sum_j \norm f_j.H^2. \norm g_j .H^2.\mid h=\sum_j  f_j g_j \Bigr\}.
\emd

One way to phrase the main result of Ferguson and Lacey \cite{sarahlacey} is that we have the equality 
\md0
H^2(\z5D)\widehat\otimes H^2(\z5D)=H^1(\z5D).
\emd
Hence, the linear functional $L$ extends to a bounded linear functional on $H^1(\z5D)$. 
By Chang--Fefferman duality, it is the case that 
 $\norm  \ZP_{-,-}\zvf.\mo B.\lesssim\norm h_\zvf..$. 

  In addition, due to the Hahn-Banach Theorem, and the inclusion $H^1\subset L^1$, we can extend 
  $L$ to all of $L^1(\z5T)$.  Hence, there is a $\zc\in L^\zI$ with  $\ZP_{-,-}\zvf=\ZP_{-,-}\zc$ and 
$\norm \zc .\zI.\lesssim{}\norm h_\zvf..$. 
\end{proof}

We remark that in the one variable case, one achieves equality in \zref e.nehari/. This is due to the fact that 
each $h\in H^1(\ZD)$ can be factored as the product of functions in $H^2$, with equality of norms.  In the present 
setting, one knows only weak factorization, that is the equality of $H^1(\z5D)$ with the injective tensor product of $H^2(\z5D)$ with itself.

To prove \zref c.compact/ we need the following lemma.  Let $S_j$ be the shift operator on $H^2(\z5D)$ associated with multiplication by $z_j$, for $j=1,2$.  

\bthm l.shift  For all compact operators $K\mid H^2_{+,+}(\z5D)\longrightarrow H_{-,-}^2(\z5D)$, we have 
$\norm KS_j^{n}S_{j'}^{m}..\longrightarrow0$, for $j,j'=1,2$. 
\ethm 

\begin{proof}  
It is enough to suppose that $j\neq j'$, for otherwise we simply have $S_{j}^{n+m}$ and only have to deal with one of the multiplication operators, and the argument we give will also work.  By symmetry we can suppose that $j=1$ and $j'=2$.  It is also enough to deal with finite rank operators since we can approximate any compact operator by finite rank operators.  We can actually check the claim for rank one operators, and only on a dense class of these operators.  So take $K$ to be defined by
$$
K(f)=\ip f,g. h\qquad\forall f\in H^{2}_{+,+}(\z5D)
$$
with $h\in H^{2}_{-,-}(\z5D)$ and $g\in H^{2}_{+,+}(\z5D)$ a polynomial of degree less than $n$ in the $z_1$ variable and less than $m$ in the $z_2$ variable.  But $(S_1^{*})^{n}(S_2^{*})^{m}g=0$, 
so  we have that $KS_1^{n}S_2^{m}=0$.
\end{proof}

\begin{proof}[Proof of  \zref c.compact/]  If $\ZP_{-,-}\zvf$ is in the $\mo B$ closure of $C(\z5T)$, 
we can choose a polynomial $\zc$,  antiholomorphic in each variable, such that 
$\norm h_{\zvf-\zc}..$ is small.  But certainly $h_\zc$ is finite rank, hence $h_\zvf$ is the norm 
limit of finite rank operators.  Hence it is compact.

Conversely, if $h_\zvf$ is compact, then for any $\ze>0$ we can choose $n$ so large that $\norm h_\zvf S_j^n..<\ze$, 
for $j=0,1,2$, where $S_0:=S_1S_2$. 
  Note  that $h_\zvf S_j^n$ is also a Hankel operator associated to the function 
\md0
\zvf_j:= {\overline z}_j^n\ZP_{-,-}{ z}_j^n\zvf,\quad j=1,2,\qquad \zvf_0:= {\overline z}_1^n {\overline z}_2^n
\ZP_{-,-}{ z}_1^n{ z}_2^n\zvf.
\emd
Thus, by  \zref t.nehari/, $\zvf_j$ has $\mo B(\z5D)$ norm at most a constant times \ze.
That is, the Hankel operator $h_\zvf$ is well approximated in operator norm by the operator associated to the  polynomial 
\md0
\sum_{-n<{}m_1,m_2<0}\widehat\zvf(m_1,m_2)z_1^{m_1}z_2^{m_2}=\ZP_{-,-}\zvf+\zvf_0-\zvf_1-\zvf_2.
\emd
Thus, we see that $\ZP_{-,-}\zvf$ is in the $\mo B$ closure of $C(\z5T)$. 
\end{proof}

\begin{proof}[Proof of \zref t.many/.]

We prove the equivalence of the compactness of the commutator $C_\zvf$ defined in \zref e.cb/ and $\zvf\in \mo V$, and prove the assertation that $\mo V^*=H^1$ in the next section. 
Central to this proof is the characteriztion of the compactness of the Hankel operators that we have already given. 
While we discussed that proof on the circle, it has an equivalent formulation on the real line.\footnote{In fact, the paper of 
Ferguson and Lacey \cite{sarahlacey} is phrased on the real line, making  certain simplifications 
 for that proof available.} 
 
 Indeed, there are four relevant Hankel operators on $L^2(\z5R)$.  They are given as maps from $H_\zve^2(\z5R) \longrightarrow
 H_{-\zve}^2(\z5R)$, where $\zve\in\{\pm,\pm\}$, and $-\zve$ is  conjugate to $\zve$.  The definition is below, with $M_\zvf$ being the 
 operator of pointwise multiplication by $\zvf$. 
 \md0
 H_{\zvf,\zve}f=\ZP_{-\zve}M_\zvf\mid H_\zve^2(\z5R) \longrightarrow
 H_{-\zve}^2(\z5R).
 \emd
 We have the fact that any  of these operators is compact iff $\ZP_{-\zve}\,\zvf\in\mo V(\z5R)$.  
 $L^2(\z5R)$ is a sum of these Hardy spaces, and  the commutator $C_\zvf$ is a sum of these four Hankel operators.  Thus, 
 $C_\zvf$ is compact iff each of the $H_{\zvf,\zve}$ are compact iff $\zvf\in \mo V(\z5R)$. 
\end{proof}

\section{The Essential Norm of Little Hankel Operators}

\begin{proof}[ Proof of \zref t.0.1/]
We first show that (ii) and (iii) are equivalent.  Suppose that
$\zvf\in\z0L^{\infty}_{-,-}+C$, then we have $\zvf=\psi+g$ with
$\psi\in\z0L^{\infty}_{-,-}$ and $g\in C$.  Then for any $f\in H^{2}(\z5D)$
we have
$$
h_{\zvf}f=\ZP_{-,-}\zvf f=\ZP_{-,-}[gf+\psi
f]=\ZP_{-,-}gf+\ZP_{-,-}\psi f=h_{g}f,
$$
with the last line following because $\psi f\in\z0L^{2}_{-,-}(\z5T)$ and
$\ZP_{-,-}(\z0L^{2}_{-,-}(\z5T))=0$.  Thus giving that (ii) implies
(iii). 

Now assume that (iii) holds, we then have a function $g\in
C(\z5T)$ such that $h_{\zvf}=h_g$.  Let $f\in H^{2}(\z5D)$, then
because $h_{\zvf}=h_g$ we have
$$
\ZP_{-,-}((\zvf-g)f)=0\qquad\forall f\in H^{2}(\z5D).
$$
Letting $f=1$ we have that $\zvf-g\in\z0L^{2}_{-,-}(\z5T)$, but we also
have that $\zvf-g\in L^{\infty}(\z5T)$, which implies that
$\zvf-g\in\z0L^{\infty}_{-,-}(\z5T)$, or $\zvf\in\z0L^{\infty}_{-,-}+C$.

Now we show that (i) and (ii) are equivalent.  But this follows
immediately from  \zref t.0.2/.  This is because $h_\zvf$ is compact if
and only if $\norm{h_\zvf}.{e}.=0$.  But if $\norm{h_{\zvf}}.e.=0$, then by
\zref t.0.2/ we have that
$\textnormal{dist}_{L^{\infty}}(\zvf,\z0L^{\infty}_{-,-}+C)=0$ and so
$\zvf\in\z0L^{\infty}_{-,-}+C$.  Conversely, if $\zvf\in\z0L^{\infty}_{-,-}+C$, then
$\textnormal{dist}_{L^{\infty}}(\zvf,\z0L^{\infty}_{-,-}+C)=0$.  But by \zref t.0.2/ we have $\norm{h_\zvf}.e.=0$, or $h_\zvf$ is compact.
\end{proof}

Our  proof of  \zref t.0.2/ is heavily influenced by the presentation of Hartman's Theorem in V.~Peller's book \cite{MR1949210}.
We will need a few
simple lemmas in the course of the proof of the theorem.

\bthm l.1.1
If $\psi\in\z0L^{\infty}_{-,-}(\z5T)$ and $\zvf\in L^{\infty}(\z5T)$ then
$h_\zvf=h_{\zvf+\psi}$.
\ethm

\begin{proof}
Let $f\in H^{2}(\z5D)$.  Then
$$h_{\zvf+\psi}f=\ZP_{-,-}(\psi+\zvf)f=\ZP_{-,-}\psi f+\ZP_{-,-}\zvf f=\ZP_{-,-}\zvf f=h_{\zvf}f,
$$
with the second to last inequality following since $\psi
f\in\z0L^{2}_{-,-}(\z5T)$ and $\ZP_{-,-}(\z0L^{2}_{-,-}(\z5T))=0$.
\end{proof}

\bthm l.1.2 
Let $\zvf\in L^{\infty}(\z5T)$.  Then
$$
\norm{h_\zvf}..\approx\inf\{\norm{\zvf-\psi}.{\infty}.:\psi\in\z0L^{\infty}_{-,-}(\z5T)\}:=\textnormal{dist}_{L^{\infty}}(\zvf,\z0L^{\infty}_{-,-}).
$$
\ethm

This is the natural extenstion to the bi-disk of the fact in one complex variable that one can approximate the norm of a Hankel operator
by the distance of its symbol from $H^{\infty}(\ZD)$.   

\begin{proof}
Clearly if $\zvf\in L^{\infty}(\z5T)$ and $\psi\in\z0L^{\infty}_{-,-}(\z5T)$
then
$$
\norm{h_\zvf}..=\norm{h_{\zvf-\psi}}..\leq\norm{\zvf-\psi}.{\infty}.,
$$
and so
$$
\norm{h_\zvf}..\leq\textnormal{dist}_{L^{\infty}}(\zvf,\z0L^{\infty}_{-,-}).
$$
On the other hand, Nehari's Theorem on the bi-disk, \zref t.nehari/ implies that
$$
\norm{\zvf-\psi}.{\infty}.\lesssim \norm{h_{\zvf-\psi}}..=\norm{h_{\zvf}}..
$$
and so $\textnormal{dist}_{L^{\infty}}(\zvf,\z0L^{\infty}_{-,-})\lesssim{}
\norm{h_{\zvf}}..$.  This proves the lemma.
\end{proof}

We are also going to need a characterization of the space $\z0L^{\infty}_{-,-}+C$.  Recall this is the space of functions $\zvf\in L^{\infty}(\z5T)$ 
that have a decomposition of the form $\psi +g$ with $\psi\in\z0L^{\infty}_{-,-}(\z5T)$ and $g\in C(\z5T)$.
Similar to the one-variable case we have the following theorem.

\bthm t.1.1
$\z0L^{\infty}_{-,-}+C$ is a closed subspace of $L^{\infty}_{-,-}(\z5T)$, and moreover 
$$\z0L^{\infty}_{-,-}+C=\textnormal{clos}_{L^{\infty}}\left(\bigcup_{n,m=0}^{\infty}\overline{z}_1^{n}\overline{z}_2^{m}\z0L^{\infty}_{-,-}(\z5T)\right).$$ 
\ethm

This result is slightly different than what one would find in one complex variable.  In one variable, the analog of this space is
$H^{\infty}+C$, which is in fact a sub-algebra of $L^{\infty}(\ZT)$.  In higher dimensions, $\z0L^{\infty}_{-,-}(\z5T)$ is not closed under 
multiplication as $H^{\infty}(\z5D)$ is, so $\z0L^{\infty}_{-,-}+C$ will not be a sub-algebra.  This is also a remnant of the fact we are 
working with little Hankel operators.  To prove this theorem, we will need one more lemma.

\bthm l.1.4
Let $C_{\z0L}(\z5T):=\z0L^{\infty}_{-,-}(\z5T)\cap C(\z5T)$ and $\zvf\in C(\z5T)$.  Then
$$
\textnormal{dist}_{L^{\infty}}(\zvf,\z0L^{\infty}_{-,-})=\textnormal{dist}_{L^{\infty}}(\zvf,C_{\z0L}).
$$
\ethm 

\begin{proof}
Since $C_{\z0L}(\z5T)\subset\z0L^{\infty}_{-,-}(\z5T)$ we trivially have that
$$
\textnormal{dist}_{L^{\infty}}(\zvf,\z0L^{\infty}_{-,-})\leq\textnormal{dist}_{L^{\infty}}(\zvf,C_{\z0L}).
$$
We are going to use the harmonic extension of functions in $L^{\infty}(\z5T)$ to the bi-disk.  The extension will also be denoted by the function element.  Finally, let $h_r(\xi)=h(r\xi)$, $\xi\in\z5T$.  Let $\psi\in\z0L^{\infty}_{-,-}(\z5T)$.  Since $\zvf\in C(\z5T)$ we have

\begin{eqnarray*}
\norm{\zvf-\psi}.{\infty}. & \geq & \lim_{r\to 1}\norm{(\zvf-\psi)_r}.{\infty}.\\
 & \geq & \lim_{r\to 1}(\norm{\zvf-\psi_r}.{\infty}.-\norm{\zvf-\zvf_r}.{\infty}.)\\
 & = & \lim_{r\to 1}\norm{\zvf-\psi_r}.{\infty}..
\end{eqnarray*}
But this then shows that $\norm{\zvf-\psi}.\infty.\geq\textnormal{dist}_{L^{\infty}}(\zvf, C_{\z0L})$ for any $\psi\in\z0L^{\infty}_{-,-}$, proving the lemma.
\end{proof}

We can now prove   \zref t.1.1/.
\begin{proof}
By \zref l.1.4/ we have that $C/C_{\z0L}$ has an isometric embedding in $L^{\infty}/ \z0L^{\infty}_{-,-}$ and can thus be considered as a closed subspace of $L^{\infty}/  \z0L^{\infty}_{-,-}$.
Let $\rho:L^{\infty}\to L^{\infty}/ \z0L^{\infty}_{-,-}$ be the natural quotient map.  Then we have that
$\z0L^{\infty}_{-,-}+C=\rho^{-1}(C/C_{\z0L})$.  This follows
from the fact that we have a subspace and are looking at the quotient map.  The proof is the same as in the one variable case.  
See \cite{MR1949210} for details.

Finally, $\z0L^{\infty}_{-,-}+C=\textnormal{clos}_{L^{\infty}}\left(\cup_{n,m=0}^{\infty}
\overline{z}_{1}^{n}\overline{z}_{2}^{m}\z0L^{\infty}_{-,-}(\z5T)\right)$.  
This follows since the continuous functions on  $\z5T$ can be uniformly approximated by polynomials in
$z_1$, $z_2$ and their conjugates.
\end{proof}

\begin{proof}[Proof of \zref t.0.2/.]
Let $K:H^{2}_{+,+}(\z5D)\to H^{2}_{-,-}(\z5D)$ be a compact operator.  Then we want to estimate $\norm{h_\zvf-K}..$ 
from below.  Let $S_j$ be multiplication by the variable $z_j$.  Using that $S_j$ is a contraction and properties of norms we have
\begin{eqnarray*}
\norm{h_\zvf-K}.. & \geq & \norm{(h_\zvf-K)S_1^{n}S_{2}^{m}}..\\
 & \geq & \norm{h_{\zvf}S_1^{n}S_2^{m}}..-\norm{KS_1^{n}S_2^{m}}..\\
  & = & \norm{h_{z_1^n z_2^m\zvf}}..-\norm{KS_1^{n}S_2^{m}}..\\
  & \gtrsim & \textnormal{dist}_{L^{\infty}}(\zvf,\overline{z}_1^{n}\overline{z}_2^{m}\z0L_{-,-}^{\infty})-\norm{KS_1^{n}S_2^{m}}..\\
  & \geq & \textnormal{dist}_{L^{\infty}}(\zvf, \z0L_{-,-}^{\infty}+C)-\norm{KS_1^{n}S_2^{m}}..
\end{eqnarray*}
We used the fact that $\norm{h_\zvf}..\approx\textnormal{dist}_{L^{\infty}}(\zvf,\z0L_{-,-}^{\infty})$ as shown in \zref l.1.2/ and the characterization of $\z0L_{-,-}^{\infty}+C$ given in
\zref t.1.1/.  Now by \zref l.shift/ we have as $n,m\to\infty$ that 
$$
\norm{KS_1^{n}S_2^{m}}..\to 0.
$$
So $\norm{h_{\zvf}-K}..\gtrsim\textnormal{dist}_{L^{\infty}}(\zvf,\z0L_{-,-}^{\infty}+C)$ for any compact operator $K$.  This then gives
$$
\norm{h_{\zvf}}.e.\gtrsim\textnormal{dist}_{L^{\infty}}(\zvf,\z0L_{-,-}^{\infty}+C).
$$

To prove the other inequality, begin by supposing that $g$ is a trigonometric polynomial.  Then $h_g$ is a compact (finite rank) operator.   So for any $g\in C(\z5T)$ the operator $h_g$ is compact.  Then we have
$$
\norm{h_{\zvf}}.e.\leq\inf_{g\in C}\norm{h_{\zvf}-h_g}..=\inf_{g\in C}\norm{h_{\zvf-g}}..\,.
$$
By \zref l.1.2/ we then have
$$
\norm{h_{\zvf}}.e.\lesssim\inf_{g\in C,\psi\in\z0L_{-,-}^{\infty}}\norm{\zvf-g-\psi}..\lesssim\textnormal{dist}_{L^{\infty}}(\zvf,\z0L_{-,-}^{\infty}+C).
$$
Combining these two estimate we have that $\norm{h_\zvf}.{e}.\approx\textnormal{dist}_{L^{\infty}}(\zvf, \z0L_{-,-}^{\infty}+C)$.
\end{proof}

\section{$\mo V $ and Carleson Measures} 
We state an equivalent form of the definition of $\mo V(\z5R)$ in terms of Carleson measures
and, in particular, in a variant of \zref e.CM/.

\bthm p.CM  Fix a choice of wavelet $w$. 
A function $b$ is in $\mo V(\z5R)$ iff any of the conditions below hold.
\begin{itemize}
\item[(i)]  $b$ is in the  closure, in $\mo B$ norm, 
 of all finite linear combinations of  
 \md0
 \{w_R\mid R\in\z0D\times \z0D\}.
 \emd
\item[(ii)] $b\in \mo B(\z5R)$, and writing $R=R_1\times R_2$ for a rectangle $R$,
\md4
\lim_{N\to\zI} \NOrm \sum_{\substack { R\in \z0D\times\z0D\\\abs{ \log\abs{R_1}}+\abs{\log \abs{R_2}}>N}}
\langle b,w_R\rangle w_R .\mo B(\z5R).=0,
\\
\lim_{N\to\zI} \NOrm \sum_{\substack { R\in \z0D\times\z0D\\
R\not\subset\{\abs x<N\}}}\langle b,w_R\rangle w_R .\mo B(\z5R).=0. 
\emd
\end{itemize}
These two  conditions are independent of the  choice of wavelet basis.
\ethm 

Set $FW(w)$ to be the linear space of finite linear combinations of $\{w_R\mid R\in\z0D\times\z0D\}$.  Our first lemma is 

\bthm l.fw  
For any two choices of $w,w'$, 
\md0
\textup{clos}_{\mo B }{\operatorname{FW}(w)}=\textup{clos}_{\mo B }{\operatorname{FW}(w')}.
\emd
\ethm 

Observe that the space $\textup{clos}_{\mo B }{\operatorname{FW}(w)}$ is invariant under dilations by factors of $2$.  And, under our assumptions on the wavelets, we have 
\md0
\sum_{R}\abs{\ip w',w_R.}<\zI.
\emd
This fact clearly implies that each wavelet $w'_R\in\textup{clos}_{\mo B }{\operatorname{FW}(w)}$, and moreover that 
the same is true of each element of ${\operatorname{FW}(w')}$.
Thus, the lemma is immediate.
This frees us to make particular choices for $w$ in different parts of our proof.  In addition, we suppress the explicit choice of wavelet in our notation.

It is a routine matter to verify that $b\in\textup{clos}_{\mo B }{\operatorname{FW}}$ iff it satisfies condition $(ii)$ of \zref p.CM/.

Let us see condition $(i)$ of the proposition, that is 
\md0
\textup{clos}_{\mo B }{C^\zI_0}=\textup{clos}_{\mo B }{FW}.
\emd
We are free to choose the wavelet to be smooth and  have compact spatial support, in which case it is clear that $FW\subset \textup{clos}_{\mo B }{C^\zI_0}$.  And so we need only argue for the reverse inclusion. 
 But it is very easy to verify that a function in $C^\zI_0$ satisfies condition (ii) of the proposition.   
 In fact, this verification depends upon the estimates below, valid for all $f\in C^\zI_0$, with constants that depend upon the choice of $f$. 
\md0
\abs{\ip f,w_R.}\seq{} 
\begin{cases}
\abs{R}^{3/2}, & \abs{R_1}+\abs{R_2}<1 
\\
\frac{\abs{R_1}}{\sqrt{\abs{R_2}}}, & \abs{R_1}<1<\abs{R_2} 
\\
 \abs{R}^{-1/2}, & \abs{R_1},\abs{R_2}>1.
\end{cases}
\emd
And so, a function in $C^\zI_0$ can be well approximated in $\mo B$ norm by finite sums of wavelets. 

\medskip 

We address the equality ${\mo V}^*=H^1$.  $H^1$ and $\mo B$ duality shows that  $H^1\subset {\mo V}^*$, and so we should show the reverse 
containment.  But duality and $\mo V=\textup{clos}_{\mo B }{FW}$ also shows that for $f\in FW$, 
\md0
\sup_{\substack{b\in \mo V \\ \norm b.\mo B.=1}}\abs{\ip f,b.}\ge{}c\norm f. H^1. .
\emd
So to conclude the identity, it would be enough to know that $\textup{clos}_{H^1 }{FW}=H^1$.  
This equality follows from one of the several equivalent definitions of $H^1$ that have been 
established by Chang and Fefferman.  
In particular, we have 
\md0
\norm f.H^1.\approx \NOrm \Bigl[ \sum_{R\in\z0R} \frac{ \abs{\ip f,w_R.}^2}{\abs R} \ind R \Bigr]^{1/2} .1..
\emd
And this equivalence proves that  $\textup{clos}_{H^1 }{FW}=H^1$.

\bigskip


\bigskip
 {\parindent=0pt 
Michael T. Lacey \hfill\break
School of Mathematics\hfill\break
Georgia Institute of Technology\hfill\break
Atlanta GA 30332\hfill\break
\smallskip
{\tt lacey@math.gatech.edu\hfill\break
\tt http://www.math.gatech.edu/\~{}lacey} \hfill\break

\bigskip 

Erin Terwilleger \hfill\break
Department of  Mathematics\hfill\break
University of Connecticut\hfill\break
Storrs, CT  06269-3009\hfill\break
\smallskip
{\tt  terwilleger@math.uconn.edu\hfill\break
\tt http://www.math.uconn.edu/\~{}terwilleger}\hfill\break

\bigskip 

Brett D. Wick \hfill\break
Department of Mathematics \hfill\break
Brown University \hfill\break 
Providence, RI 02912\hfill\break 
\smallskip 
{\tt bwick@math.brown.edu\hfill\break
\tt http://www.math.brown.edu/\~{}bwick}


\section*{References} 

\begin{biblist}


\bib{bourdaud}{article}{
    author={Bourdaud, G{\'e}rard},
     title={Remarques sur certains sous-espaces de ${\rm BMO}(\mathbb R\sp n)$
            et de ${\rm bmo}(\mathbb R\sp n)$},
  language={French},
   journal={Ann. Inst. Fourier (Grenoble)},
    volume={52},
      date={2002},
    number={4},
     pages={1187\ndash 1218},
      issn={0373-0956},
    review={1 927 078},
}






\bib{cf1}{article}{
    author={Chang, Sun-Yung A.},
    author={Fefferman, Robert},
     title={Some recent developments in Fourier analysis and $H\sp p$-theory
            on product domains},
   journal={Bull. Amer. Math. Soc. (N.S.)},
    volume={12},
      date={1985},
    number={1},
     pages={1\ndash 43},
      issn={0273-0979},
    review={MR 86g:42038},
}
 
\bib{cf2}{article}{
    author={Chang, Sun-Yung A.},
    author={Fefferman, Robert},
     title={A continuous version of duality of $H\sp{1}$ with BMO on the
            bidisc},
   journal={Ann. of Math. (2)},
    volume={112},
      date={1980},
    number={1},
     pages={179\ndash 201},
      issn={0003-486X},
    review={MR 82a:32009},
}


\bib{MR54:843}{article}{
    author={Coifman, R. R.},
    author={Rochberg, R.},
    author={Weiss, Guido},
     title={Factorization theorems for Hardy spaces in several variables},
   journal={Ann. of Math. (2)},
    volume={103},
      date={1976},
    number={3},
     pages={611\ndash 635},
    review={MR 54 \#843},
}
\bib{sarahlacey}{article}{
    author={Ferguson, Sarah H.},
    author={Lacey, Michael T.},
     title={A characterization of product BMO by commutators},
   journal={Acta Math.},
    volume={189},
      date={2002},
    number={2},
     pages={143\ndash 160},
      issn={0001-5962},
    review={1 961 195},
}

\bib{fergusonsadosky}{article}{
    author={Ferguson, Sarah H.},
    author={Sadosky, Cora},
     title={Characterizations of bounded mean oscillation on the polydisk in
            terms of Hankel operators and Carleson measures},
   journal={J. Anal. Math.},
    volume={81},
      date={2000},
     pages={239\ndash 267},
      issn={0021-7670},
    review={MR 2001h:47040},
}


\bib{nehari}{article}{
    author={Nehari, Zeev},
     title={On bounded bilinear forms},
   journal={Ann. of Math. (2)},
    volume={65},
      date={1957},
     pages={153\ndash 162},
    review={MR 18,633f},
}

\bib{MR1949210}{book}{
    author={Peller, Vladimir V.},
     title={Hankel operators and their applications},
    series={Springer Monographs in Mathematics},
 publisher={Springer-Verlag},
     place={New York},
      date={2003},
     pages={xvi+784},
      isbn={0-387-95548-8},
    review={1 949 210},
}

\bib{sarason}{article}{
    author={Sarason, Donald},
     title={Functions of vanishing mean oscillation},
   journal={Trans. Amer. Math. Soc.},
    volume={207},
      date={1975},
     pages={391\ndash 405},
    review={MR 51 \#13690},
}
 



 
\end{biblist}
\end{document}